# Vertical Trajectory Optimization to Minimize Environmental Impact in the Presence of Wind

Sang Gyun Park[*],
John-Paul Clarke[†]
*Georgia Institute of Technology, Atlanta, GA 30332*

The vertical trajectory optimization for the en route descent phase is studied in the presence of both along track and cross winds, which are both modeled as functions of altitude. The flight range covers some portion of a cruise segment and ends at a meter fix. The descent trajectory is assumed as a flight idle descent. The problem is formulated as an optimal control problem with both mixed state path constraints and pure state constraints. For minimizing environmental impacts, we optimize descent trajectory with respect to two cost functionals: fuel and emissions cost. We analyze both singular arc and boundary arc using the necessary conditions of the optimality. From the analysis result, we propose the optimal trajectory generation method, in which the optimal trajectory is generated by the backward and forward integration. The trajectory from the proposed algorithm is compared to the numerical optimal solution of the original optimal control problem. The result shows that proposed algorithm generates the same solution as the numerical optimal solution. The wind speed, wind shear, and cross wind effects on the optimal trajectory are analyzed with two aircraft type, Boeing 737-500 and Boeing 767-400.

## Nomenclature

| | | |
|---|---|---|
| $D$ | = | Drag force |
| $h$ | = | Altitude |
| $K_{cr}$ | = | Cruise cost |
| $K_{des}$ | = | Idle descent cost |
| $L$ | = | Lift force |
| $M$ | = | Mach number |
| $m$ | = | Aircraft mass |
| $T$ | = | Thrust force |
| $V_{CAS}$ | = | Calibrated airspeed |
| $V_G$ | = | Ground speed |
| $V_T$ | = | True airspeed |
| $W_h$ | = | Along track component of the wind velocity |
| $W_c$ | = | Cross track component of the wind velocity |
| $x_s$ | = | Along track distance |
| $()^W$ | = | Vector expressed in relative wind frame |
| $\gamma$ | = | Aerodynamic flight path angle |
| $\psi$ | = | Heading angle |
| $\psi_w$ | = | Angle between airspeed and ground speed |

[*]Ph.D. , currently working at Optimal Synthesis Inc.,
[†]Professor, School of Aerospace Engineering

## I. Introduction

THE Next Generation Air Transportation System (NextGen) is being developed and implemented to enhance the efficiency and safety of air transportation and reduce the environmental impact due to noise and aircraft engine gaseous emissions. The Continuous Descent Arrival (CDA) or Optimized Profile Descent (OPD) is a promising approach to reduce the environmental impact of aircraft operations[1, 2, 3]. During CDA\OPD, aircraft descend from the cruise altitude to the Final Approach Fix (FAF) at or near idle thrust without level segments at low altitude, thereby minimizing the need for high thrust levels to remain at a constant altitude, and reducing the environmental impact[4].

In order to maximize the environmental benefits of CDA, Park and Clarke [5, 6] formulated the vertical trajectory optimization problem as a multi-phase optimal control problem and showed that the possible trajectory variation with respect to the various performance indices occurs in en route descent phase, which starts during the latter stages of the cruise segment and ends at the meter fix. Thus, the focus of this paper is the vertical trajectory generation problem in the en route descent phase.

Optimal control techniques have been widely used to solve trajectory generation problems. In 1980s, most approaches used energy state approximation in which the energy state was defined by the combination of height and speed to reduce the problem size[7, 8, 9, 10, 11]. Erzberger and Lee[7] proposed an optimal trajectory generation method with specified range. They used Direct Operational Cost (DOC) combining time cost and fuel cost as a performance index. Sorenson and Waters[8] addressed the fuel optimal trajectory problem with fixed arrival time. Available time delay analysis was performed by setting negative time cost. Burrows[9] converted the fixed arrival time problem to a DOC optimal trajectory problem with free final time. Using this approach, the fuel optimal control problem is equivalent to finding a time cost for which a corresponding free final time DOC optimal trajectory arrives at the assigned time. Chakravarty[10] used the singular perturbation method as well as an energy state approximation and investigated the sensitivity of the fuel optimal trajectories to wind. Shultz[11] studied the three dimensional minimum time problem with fixed initial and final points.

Energy state approximation is beneficial because it enables quicker computation and thus lower solution times. However, since the model is simplified, and some constraints in the real procedure are not considered, trajectories generated by these approaches are typically not accurate enough to apply directly in real air traffic situations.

From early 2000s, several studies for trajectory optimization that minimize environmental impact have been conducted. Visser and Wijnen [12] addressed noise abatement trajectory optimization problem. They formulated this problem as a multi-phase optimal control problem and solved by using direct numerical method. Wu and Zhao [13] assumed descent trajectory with several segments and converted optimal control problem to parameter optimization problems with fuel and emission costs. Few papers used CDA trajectory structure assuming idle thrust descent [14, 5, 15]. Franco et al.[14] solved maximum descent range problem by singular arc analysis in the presence of the along track wind. However, they did not consider the cross wind and the flight envelope which is modeled as a pure state inequality constraint. Zhao and Tsiotras [16] proposed optimal speed profile generation method for minimizing energy with given path using singular arc analysis. However, they did not consider the wind effect.

This paper considers the optimal trajectory that minimizes environmental impact such as fuel burn and emissions. The flight range considered in this paper is also fixed from initial point in the cruise segment to the meter fix determined by Standard Terminal Arrival Route (STAR) procedure. In the previous work from same authors, optimal control framework for vertical trajectory optimization with various performance indices was developed[5], and fuel optimal trajectory generation method using Flight Management System (FMS) Vertical Navigation (VNAV) function was developed and solved in the hybrid system framework[17]. In this paper, we analyze the optimal solutions with respect to the fuel and emission performance indices by using the necessary conditions of the optimal control problem with state inequality path constraints.

The remaining part of the paper is structured as follows: The CDA trajectory structure and optimal control problem are presented in section II. The analysis of the optimal solution structure including the singular arc and boundary arc are presented in section III. The fast trajectory generation method is presented in section IV. Numerical examples with various wind conditions are presented in section V. The conclusions of this study are presented in section VI.

## II. Optimal Control Problem Formulation

### A. CDA trajectory structure

The CDA trajectory structure is shown in Fig. 1. As shown in the figure, the CDA procedure begins towards the end of the cruise segment and ends at either the FAF or the Initial Approach Fix (IAF). The traffic flow is metered at the entry points of the terminal area, where speed and altitude constraints are typically prescribed. Therefore, the CDA trajectory can be divided into two parts: an en route descent phase and a terminal area descent and approach. Only the en route descent is considered in this paper. It consists of two segments: cruise and descent. Typically, the lateral path is defined via a STAR that provides the waypoints for the Lateral Navigation (LNAV) function in the FMS [18]. Thus, only the vertical path is free to be optimized.

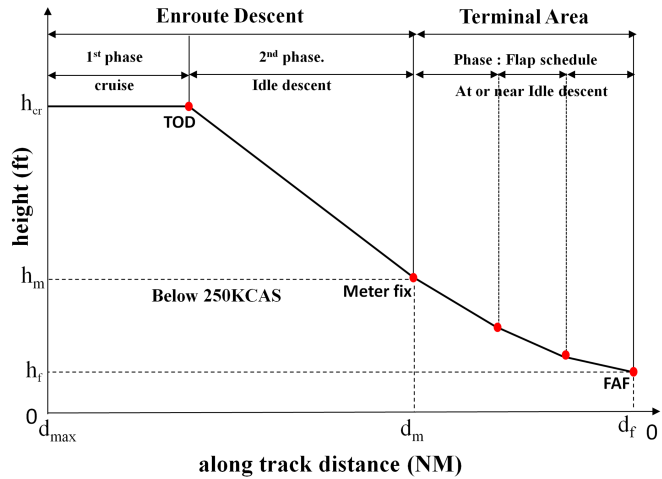

Figure 1. CDA trajectory structure

### B. Flight Dynamic Model

The relation between airspeed, cross wind, and along track wind is defined in Fig. 2. The angles $\psi$ and $\psi_w$ denote the heading angle and the angle between airspeed and ground speed, respectively.

As previously stated, only vertical motion need be considered for trajectory optimization because the lateral path of an arrival procedure is typically prescribed. In addition, if the flight path angle is bounded by small values between $-6°$ and $0°$, $\sin\gamma$ and $\cos\gamma$ can be approximated as $\gamma$ and 1, respectively, and the equations of motion in the vertical plane in [5] can be simplified as follows:

$$\dot{V}_T = \frac{1}{m}(T - D) - g\gamma - V_T\gamma\left(c(V_T, h)\frac{dW_h}{dh} + s(V_T, h)\frac{dW_c}{dh}\right) \quad (1)$$

$$\dot{x}_s = c(V_T, h)V_T + W_h \quad (2)$$

$$\dot{h} = V_T\gamma \quad (3)$$

where, $V_T$ is the true airspeed of aircraft, $x_s$ is the along track distance from runway threshold; $h$ is the altitude; $\gamma$ is the aerodynamic flight path angle; $W_h$ and $W_c$ are the along track wind speed and cross wind speed, respectively; $dW_h/dh$ and $dW_c/dh$ are the wind shear terms, $c(V_T, h) = \cos \psi_w$, and $s(V_T, h) = \sin \psi_w$. Both the along track wind and the cross wind are assumed to be functions of altitude. The aircraft mass $m$ is assumed to be a constant during the entire en route descent phase. The lift force is assumed to be equal to the gravity force. The simplified lift force $L$ with given assumption and drag force $D$ are as follows:

$$L = \frac{1}{2}\rho(h)V_T^2 S C_L = mg \qquad (4)$$

$$D(h, V_T) = \frac{1}{2}\rho(h)V_T^2 S C_D(V_T, h, C_L) \qquad (5)$$

where $\rho$ is the air density which is the function of $h$, and $C_L$ and $C_D$ are the lift and drag coefficients, respectively. The International Standard Atmosphere (ISA) model is used. From Eq. (4), $C_L = C_L(V_T, h)$. Hence, the drag in Eq. (5) is a function of $V_T$ and $h$ also.

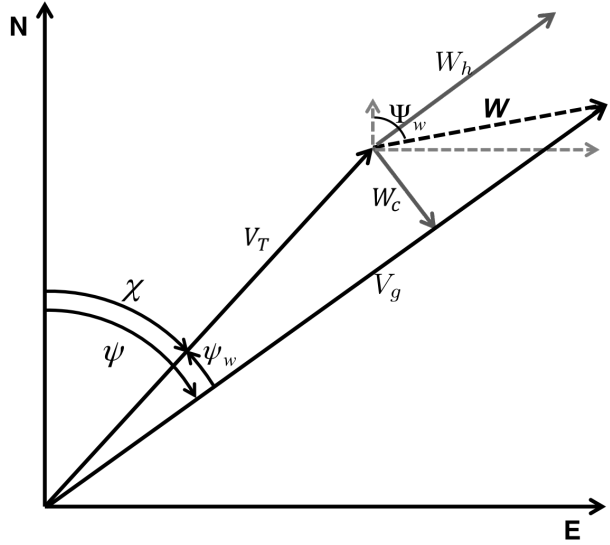

**Figure 2. Relation between airspeed, ground speed, and wind vector**

During the descent segment, thrust $T$ is set to idle as shown in Fig. 1. Generally, idle thrust is modeled as a function of $V_T$ and $h$ using data from aircraft performance models such as the performance engineering manual provided to airlines by the manufacturers and the Base of Aircraft Data (BADA) from Eurocontrol[19]. Therefore, $\gamma$ is the only control input during the idle descent segment.

### C. Environmental Cost Indices

In this study, fuel burn cost and emissions cost are used as environmental cost indices. The cost functional is defined as a summation of the cost of the cruise segment and that of the idle descent segment. From the wind speed assumption, the ground speed in the cruise phase is a constant. Therefore, the cost of the cruise phase can be expressed as a Mayer cost term. For this reason, the optimal control problem for en route descent part of CDA is formulated as a single phase optimal control problem with a cost functional of the form

$$J = K_{cr}(x_s(t_0) - d_{max}) + \int_{t_0}^{t_f} K_{des}(V_T, h) \, dt \qquad (6)$$

where $K_{cr}$ is the cruise cost coefficient which is the cost per distance, and $K_{des}$ is the idle descent cost coefficient which is the cost per time; $d_{max}$ is the along track distance from runway to the initial point;

$x_s(t_0)$ is the Top Of Descent (TOD) point in Fig. 1. The fuel cost and emissions cost can be expressed as follows:

- Fuel burn cost

$$K_{cr} = \frac{\dot{f}_{cr}}{V_{cr}} \tag{7}$$

$$K_{des}(V_T, h) = \dot{f}_{idle}(V_T, h)$$

- Emissions cost

$$K_{cr} = \frac{EI_X \dot{f}_{cr}}{V_{cr}} \tag{8}$$

$$K_{des}(V_T, h) = EI_X(V_T, h)\dot{f}_{idle}(V_T, h)$$

where $f_{cr}$ and $f_{idle}$ denote fuel flow rate in the cruise segment and the idle descent segment, respectively; $EI_X$ denotes emission index, and subscript $X$ denotes emission gas types: $CO$, $HC$, and $NO_x$.

The emissions index table can be obtained from the International Civil Aviation Organization (ICAO) engine exhaust emissions databank[20], which is subsequently interpolated using Boeing method2 (BM2), an $EI$ correction method considering airspeed and atmosphere conditions such as temperature, pressure, and humidity[21]. Since the idle fuel flow rate is a function of Mach number and $h$, $EI$ obtained by BM2 is also a function of $V_T$ and $h$. Therefore, $K_{des}$ in Eq. (6) can be expressed as $K_{des}(V_T, h)$ in both fuel and emission cost cases.

### D. Constraints for CDA

We consider several path constraints for flight envelope protection and passenger comfort. Bounds are placed on the Calibrated airspeed (CAS) and the Mach number for flight envelope protection. CAS and Mach are functions of $V_T$ and $h$ as below:

$$V_{CAS}(V_T, h) = \left[ 7R(\Theta_0)_{ISA} \left\{ \left(1 + P_\delta(h) \left[ \left(1 + \frac{V_T^2}{7R\Theta(h)}\right)^{3.5} - 1 \right] \right)^{\frac{2}{7}} - 1 \right\} \right]^{0.5} \tag{9}$$

$$M(V_T, h) = V_T / \sqrt{1.4R\Theta(h)} \tag{10}$$

where $P_\delta$ is a pressure ratio between current altitude and mean sea level (MSL); $R$ is a gas constant; $(\Theta_0)_{ISA}$ and $\Theta$ denote the standard atmospheric temperature at MSL and the temperature at the current altitude, respectively. Hence, Mach/CAS bounds are pure state inequality constraints.

For passenger comfort, we add a Rate of Descent (ROD) bound and a flight path angle bound. Boundary conditions are given in Eq. (15) and (16). CAS and altitude constraints are usually given at meter fixes per the published STAR chart. CAS constraint with an altitude constraint can be converted to the fixed $V_T$ condition at the meter fix.

- Flight Envelope Protection

$$V_{min,CAS} \leq V_{CAS} \leq V_{max,CAS} \tag{11}$$
$$M_{min} \leq M \leq M_{max} \tag{12}$$

- Passenger Comfort

$$ROD_{min} \leq -\frac{dh}{dt} \leq ROD_{max} \tag{13}$$
$$\gamma_{min} \leq \gamma \leq \gamma_{max} \tag{14}$$

- Boundary conditions

$$V_T(t_0) = V_{T0}, \quad h(t_0) = h_0, \quad x_s(t_0) = \text{free}, \quad t_0 = 0 \quad (15)$$
$$V_T(t_f) = V_{Tf}, \quad h(t_f) = h_f, \quad x_s(t_f) = s_f, \quad t_f = \text{free} \quad (16)$$

## III. Analysis of the Optimal Trajectory

### A. Equivalent Optimal Control Problem

Consider the following optimal control problem:

**Problem 1.**

$$\min_{\gamma} \quad J = \int_{t_0}^{t_f} -K_{cr}(cV_T + W_h) + K_{des}(V_T, h) \, dt$$

s.t. (17)

dynamic constraints: Eq. (1) $\sim$ (3)

path constraints : Eq. (11) $\sim$ (14)

boundary conditions : Eq. (15) and (16)

The following lemma shows that Problem 1 is equivalent to the original problem with Eq. (6).

**Lemma 1.** *Let $x_s^*(t_0)$ be the optimal TOD of Problem 1. If $x_s^*(t_0) \leq d_{max}$, Problem 1 is equivalent to the original optimal problem with performance index in Eq. (6).*

*Proof.* From Eq. (2), $x_s(t_0)$ can be expressed as

$$x_s(t_0) = x_s(t_f) - \int_{t_0}^{t_f} cV_T + W_h \, dt. \quad (18)$$

Substituting Eq. (18) into Eq. (6), then

$$J = K_{cr}(x_s(t_f) - d_{max}) + \int_{t_0}^{t_f} -K_{cr}(cV_T + W_h) + K_{des}(V_T, h) \, dt \quad (19)$$

Since $K_{cr}(x_s(t_f) - d_{max})$ is a constant from the terminal condition in Eq. (16), Problem 1 is equivalent to the original problem with Eq. (6). $\square$

Note that in the last part of this paper, we use Problem 1 to derive the necessary conditions of the original problem with cost in (6) from lemma 1.

### B. Necessary Conditions

The Hamiltonian for Problem 1 is defined as follows:

$$H = -K_{cr}(cV_T + W_h) + K_{des}(V_T, h) + \lambda_V \left( \frac{T - D}{m} - g\gamma - V_T \gamma \mathbf{W}_{h,\chi} \right) + \lambda_x(cV_T + W_h) + \lambda_h V_T \gamma \quad (20)$$

where $\lambda_V, \lambda_x,$ and $\lambda_h$ are the costate variables corresponding to the state variables $V_T$, $x_s$ and $h$, respectively. In Eq. (20), $\mathbf{W}_{h,\chi} = c\frac{dW_h}{dh} + s\frac{dW_c}{dh}$. In this paper, we define the Lagrangian using the direct adjoining approach in [22]. The Lagrangian $L$ including path constraints is given by

$$L = H + \mu^T C(V_T, h, \gamma) + \eta^T S(V_T, h) \quad (21)$$

where the inequality constraints $C(V_T, h, \gamma) : \mathbb{R}^3 \to \mathbb{R}^q$ and $S(V_T, h) : \mathbb{R}^2 \to \mathbb{R}^s$ represent the mixed path constraints and the pure state inequality in Problem 1, respectively; $\mu \in \mathbb{R}^q$ and $\eta \in \mathbb{R}^s$ are the Lagrange multiplier for mixed and pure state inequality constraints, respectively. From the constraints in Eqs. (11) $\sim$ (14), $C(V_T, h, \gamma)$ and $S(V_T, h)$ have four components as follows:

$$C = \begin{bmatrix} -\dot{h}(V_T, \gamma) - ROD_{max} \\ ROD_{min} + \dot{h}(V_T, \gamma) \\ \gamma - \gamma_{max} \\ \gamma_{min} - \gamma \end{bmatrix} \leq 0, \quad S = \begin{bmatrix} V_{CAS}(V_T, h) - V_{max,CAS} \\ V_{min,CAS} - V_{CAS}(V_T, h) \\ M(V_T, h) - M_{max} \\ M_{min} - M(V_T, h) \end{bmatrix} \leq 0 \quad (22)$$

The following equations are the necessary conditions for the optimality of the Problem 1[23].

From the Euler-Lagrange equations,

$$\dot{\lambda}_V = -\frac{\partial L}{\partial V_T} = (K_{cr} - \lambda_x)\left(c + \frac{\partial c}{\partial V_T} V_T\right) - \frac{\partial K_{des}}{\partial V_T} + \lambda_V(\frac{1}{m}\frac{\partial \tilde{D}}{\partial V_T} + \mathbf{W}_{h,\chi}\gamma + \frac{\partial \mathbf{W}_{h,\chi}}{\partial V_T} V_T \gamma) - \lambda_h \gamma - \mu^T \frac{\partial C}{\partial V_T} - \eta^T \frac{\partial S}{\partial V_T} \quad (23)$$

$$\dot{\lambda}_x = -\frac{\partial L}{\partial x_s} = 0 \quad (24)$$

$$\dot{\lambda}_h = -\frac{\partial L}{\partial h} = (K_{cr} - \lambda_x)\left(\frac{\partial c}{\partial h} V_T + \frac{dW_h}{dh}\right) - \frac{\partial K_{des}}{\partial h} + \lambda_V(\frac{1}{m}\frac{\partial \tilde{D}}{\partial h} + V_T \gamma \frac{\partial \mathbf{W}_{h,\chi}}{\partial h}) - \mu^T \frac{\partial C}{\partial h} - \eta^T \frac{\partial S}{\partial h} \quad (25)$$

where $\tilde{D} = D - T$, which means the net drag force.

From the Karush-Kuhn-Tucker (KKT) conditions,

$$\begin{aligned} \mu \geq 0, \quad & \mu^T C(V_T, h, \gamma) = 0, \\ \eta \geq 0, \quad & \eta^T S(V_T, h) = 0. \end{aligned} \quad (26)$$

The optimal control is determined by

$$\gamma^o = \mathrm{argmin}_{\gamma \in \Omega(V_T, h)} H \quad (27)$$

$$L_\gamma = H_\gamma + \mu^T C_\gamma = 0 \quad (28)$$

where $\Omega(V_T, h) = \{\gamma \mid C(V_T, h, \gamma) \leq 0\}$ is the admissible control set, and it depends on the state variable $V_T$ and $h$. Since $x_s(t_0)$ and final time $t_f$ are free, the following transversality conditions hold:

$$\lambda_x(t_0) = 0 \quad (29)$$

$$H(t_f) = 0 \quad (30)$$

Since state inequalities are independent to time and $x_s$, for any time $\tau$ on the state boundary arc where $S(V_T, h) = 0$, the following jump conditions hold[23]:

$$\lambda_x(\tau^-) = \lambda_x(\tau^+) \quad (31)$$

$$H(\tau^-) = H(\tau^+) \quad (32)$$

Therefore $\lambda_x$ and $H$ are continuous function with respect to time. Combining with Eq. (24), this condition leads to

$$\lambda_x(t) = 0 \quad \text{for} \quad t \in [t_0, t_f]. \quad (33)$$

Furthermore, since the Hamiltonian in Eq. (21) is not an explicit function of time and the cost index in Eq. (17) is in a Lagrange form, the Hamiltonian is a constant along the optimal trajectory. Hence, from Eq. (30),

$$H(t) = 0 \quad \text{for} \quad t \in [t_0, t_f] \quad (34)$$

However, $\lambda_V$ and $\lambda_h$ can be discontinuous for any time $\tau$ in the boundary arc interval and for any contact time $\tau$ by the following jump conditions[23]:

$$\lambda_V(\tau^-) = \lambda_V(\tau^+) - \nu(\tau)^T \frac{\partial S}{\partial V_T} \tag{35}$$

$$\lambda_h(\tau^-) = \lambda_h(\tau^+) - \nu(\tau)^T \frac{\partial S}{\partial h} \tag{36}$$

$$\nu(\tau) \geq 0, \quad \nu(\tau) S(\tau) = 0. \tag{37}$$

To solve this problem with the necessary conditions in Eq. (23)$\sim$ (34), both the interior arc on which $S(V_T, h) < 0$ and the boundary arc on which at least one of the components of $S(V_T, h)$ is zero must be considered[24].

## C.  Singular Arc

First, assume that all the state constraints $S(V_T, h)$ are not active. Then $\eta$ is set to zero. By the minimum principle in Eq. (27), the optimal control input $\gamma$ is determined as

$$\gamma^o = \begin{cases} \gamma_{\max} & \text{if } H_\gamma < 0 \\ \tilde{\gamma} & \text{if } H_\gamma = 0 \\ \gamma_{\min} & \text{if } H_\gamma > 0 \end{cases} \tag{38}$$

where $\tilde{\gamma}$ is the singular control, and $\gamma_{\max}$ and $\gamma_{\min}$ are determined by $\Omega(V_T, h)$. The switching function $H_\gamma$ in Eq. (38) is given by

$$H_\gamma = -\lambda_V(g + V_T \mathbf{W}_{h,\chi}) + \lambda_h V_T. \tag{39}$$

For the analysis of the singular arc, assume that $H_\gamma$ is zero during a finite time interval. Then the singular control and the singular arc are obtained by the time derivatives of $H_\gamma$. From $H_\gamma = 0$, we can obtain the following relation:

$$(g + V_T \mathbf{W}_{h,\chi}) = \frac{\lambda_h}{\lambda_V} V_T \tag{40}$$

By differentiating Eq. (39), and substituting for the time derivatives of the costate and state variables and Eq. (40), $\dot{H}_\gamma$ is expressed as

$$\begin{aligned}
\dot{H}_\gamma = K_{cr} \left[ -\left(c + \frac{\partial c}{\partial V_T} V_T\right) \frac{\lambda_h}{\lambda_V} V_T + \frac{\partial c}{\partial h} V_T^2 + V_T \mathbf{W}_{h,\chi} \right] + \frac{\partial K_{des}}{\partial V_T} \frac{\lambda_h}{\lambda_V} V_T - \frac{\partial K_{des}}{\partial h} V_T \\
+ \frac{\lambda_V}{m} \left[ \tilde{D} \left( \mathbf{W}_{h,\chi} + \frac{\partial \mathbf{W}_{h,\chi}}{\partial V_T} \right) + \frac{\partial \tilde{D}}{\partial h} V_T \right] - \frac{\lambda_h}{m} \left( \frac{\partial \tilde{D}}{\partial V_T} V_T + \tilde{D} \right) = 0.
\end{aligned} \tag{41}$$

From Eq. (34) and the singular control assumption $H_\gamma = 0$, the remaining term in Eq. (20) is zero on the singular arc also. Therefore,

$$H_0 = -K_{cr}(cV_T + W_h) + K_{des} - \lambda_V \frac{\tilde{D}}{m} = 0 \tag{42}$$

By combining Eq. (42) and Eq. (40),

$$\begin{aligned}
\frac{\lambda_V}{m} &= \frac{1}{\tilde{D}}(-K_{cr}(cV_T + W_h) + K_{des}) \\
\frac{\lambda_h}{m} &= \frac{\lambda_V}{m} \frac{1}{V_T}(g + V_T \mathbf{W}_{h,\chi}).
\end{aligned} \tag{43}$$

By substituting Eq. (43) into Eq. (41), the final formulation of the singular arc $\Gamma_s(V_T, h)$ is obtained as

$$\Gamma_s(V_T, h) = [K_{cr}(cV_T + W_h) - K_{des}(V_T, h)] \left[ V_T \frac{\partial \tilde{D}}{\partial h} - \frac{\partial \tilde{D}}{\partial V_T}(g + V_T \mathbf{W}_{h,\chi}) + \tilde{D} \frac{\partial \mathbf{W}_{h,\chi}}{\partial V_T} \right]$$
$$- \tilde{D} \left[ K_{cr} \left( cg + \frac{W_h}{V_T} g - \left( c + \frac{\partial c}{\partial V_T} V_T \right)(g + V_T \mathbf{W}_{h,\chi}) + \frac{\partial c}{\partial h} V_T^2 + V_T \mathbf{W}_{h,\chi} \right) \right]$$
$$- \tilde{D} \left[ -\frac{K_{des}}{V_T} g + \frac{\partial K_{des}}{\partial V_T}(g + V_T \mathbf{W}_{h,\chi}) - \frac{\partial K_{des}}{\partial h} V_T \right]$$
$$= 0$$
(44)

Since the singular arc $\Gamma_s(V_T, h)$ in Eq. (44) is not an explicit function of control $\gamma$, second time derivative $\ddot{H}_\gamma$ is needed to obtain the analytic formulation of optimal singular control $\tilde{\gamma}$. Furthermore, the following Generalized Legendre-Clebsch (GLC) condition should hold if the singular arc is a part of the optimal trajectory.

$$\frac{\partial}{\partial \gamma} \left[ \frac{d^2}{dt^2} H_\gamma \right] = \frac{\partial}{\partial \gamma} \left[ \frac{d}{dt} \dot{H}_\gamma \right]$$
$$= \frac{\partial}{\partial \gamma} \left[ \frac{\partial \dot{H}_\gamma}{\partial V_T} \left( -\frac{\tilde{D}}{m} - (g + V_T \mathbf{W}_{h,\chi})\gamma \right) + \frac{\partial \dot{H}_\gamma}{\partial h}(V_T \gamma) \right] \quad (45)$$
$$= -\frac{\partial \dot{H}_\gamma}{\partial V_T}(g + V_T \mathbf{W}_{h,\chi}) + \frac{\partial \dot{H}_\gamma}{\partial h} V_T \leq 0$$

Here, the partial derivative terms of $\dot{H}_\gamma$ can be calculated numerically when Eq. (44) is solved. Hence, this condition can be checked on the singular arc. All numerical examples in this paper satisfy this condition. Furthermore, if the singular control holds $\gamma_{\min} < \tilde{\gamma} < \gamma_{\max}$, the singular arc control satisfies the following junction theorem also.

**Theorem 1** (Junction Thm in [25]). *Let $t_c$ be a point at which singular and nonsingular sub arcs of an optimal control $u$ are joined, and let $q$ be the order of the singular arc. Suppose the strengthened GLC condition is satisfied at $t_c$, and assume that the control is piecewise analytic in the neighborhood of $t_c$. Let $u^{(r)}$ be the lowest order derivative of u which is discontinuous at $t_c$. Then $q + r$ is an odd integer.*

In this problem, $q$ is one because $H_\gamma^{(2)} = \ddot{H}_\gamma$ has the explicit formula with respect to $\gamma$. The control $\gamma$ in the nonsingular arc has either $\gamma_{\max}$ or $\gamma_{\min}$, and hence $\gamma$ is discontinuous at the junction point, which means $r = 0$. Therefore, $q + r = 1$ in this problem.

The singular arcs of B737-500 aircraft with various wind conditions are shown in Fig. 3. The right sides of singular arcs have negative values of $\Gamma_s(V_T, h)$ for both the fuel cost and the $NO_x$ cost cases. The singular controls of B737-500 with various wind conditions are shown in Fig. 4. As may be seen, discontinuity assumption at the junction point is true for any point on the singular arc.

### D. Boundary Arc

In this subsection, we analyze the boundary arc, on which at least one of $S(V_T, h)$ is zero. To do this, we assume that there exists a nonzero time interval $[t_1, t_2]$ in which the optimal trajectory is on the active state inequality. We also assume that the mixed state inequality $C(V_T, h, \gamma)$ is not active on the boundary arc. Therefore, $\mu(t) = 0, t \in [t_1, t_2]$. This assumption is reasonable because the flight path angle during the idle descent flight with maximum/minimum constant MACH/CAS is within the bound of the flight path angle limit.

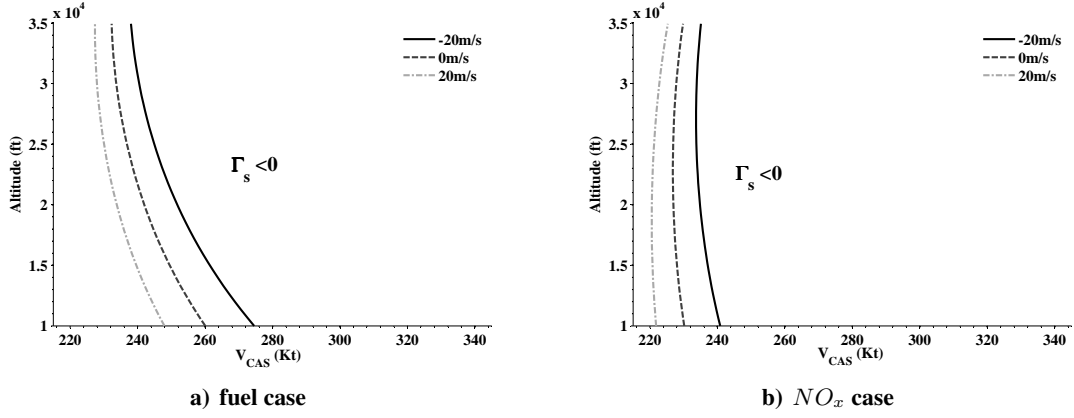

a) fuel case    b) $NO_x$ case

Figure 3.  B737-500 Singular arcs in various wind conditions

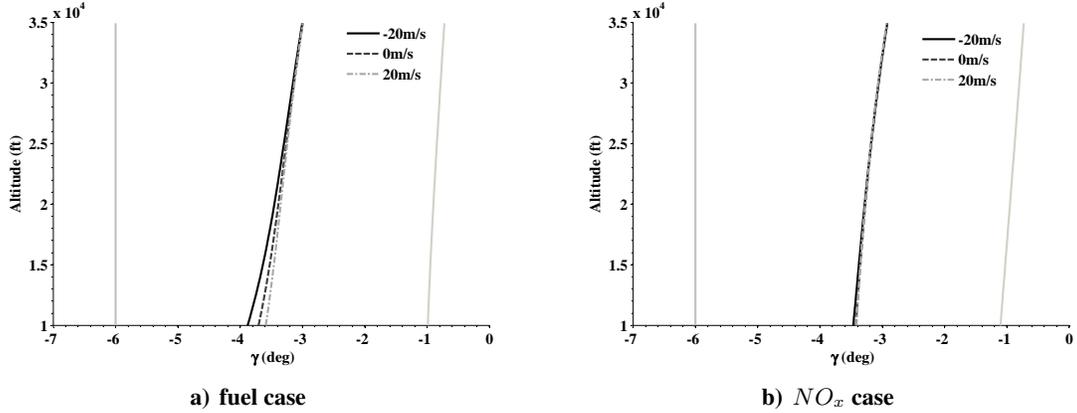

a) fuel case    b) $NO_x$ case

Figure 4.  B737-500 singular controls in various wind conditions

From the necessary conditions for optimality and $\mu = 0$, the following equations must be satisfied for the boundary arc interval:

$$\begin{aligned}
L_\gamma(V_T, h, \gamma, \lambda_V, \lambda_h) = H_\gamma &= 0 \quad \text{in Eq. (39)} \\
H_0(V_T, h, \lambda_V, \lambda_h) &= 0 \quad \text{in Eq. (42)} \\
S_a(V_T, h) &= 0 \\
\dot{S}_a(V_T, t, \gamma) &= 0 \\
\dot{\lambda}_V &= -\frac{\partial H}{\partial V_T} - \eta_a \frac{\partial S_a}{\partial V_T} \\
\dot{\lambda}_h &= -\frac{\partial H}{\partial h} - \eta_a \frac{\partial S_a}{\partial h} \\
\eta_a(t) \geq 0, \quad t &\in [t_1, t_2]
\end{aligned} \qquad (46)$$

The active boundary arc $S_a = 0$ gives the relation between $V_T$ and $h$. Since $S_a(V_T, h)$ is the first order state inequality in this problem, which means that $\gamma$ related terms appear in the first time derivative $\dot{S}_a$, $\gamma$ on the boundary arc can be obtained from the $\dot{S}_a = 0$, and it has an explicit formulation as

$$\gamma = -\left[\frac{\partial S_a}{\partial V_T}(g + V_T \mathbf{W}_{h,\chi}) - \frac{\partial S_a}{\partial h} V_T\right]^{-1} \frac{\partial S_a}{\partial V_T} \frac{\tilde{D}}{m}. \qquad (47)$$

The Lagrange multiplier of the active state inequality $\eta_a$ is calculated by substituting the time derivative of $\lambda_V$ obtained by differentiating $\lambda_V(V_T, h)$ into the adjoint equations of motion in Eq (46)

$$\eta_a = -(\dot{\lambda}_V + \frac{\partial H}{\partial V_T})\frac{\partial S_a}{\partial V_T}^{-1}$$
$$= -\frac{m\gamma}{\tilde{D}^2}\Gamma_s(V_T, h)\frac{\partial S_a}{\partial V_T}^{-1} \geq 0 \quad (48)$$

where $\Gamma_s(V_T, h)$ is in Eq. (44). However, $\Gamma_s(V_T, h)$ can have a nonzero value on the boundary arc since $\Gamma_s(V_T, h) = 0$ only on the interior singular arc.

In Eq. (48), $\Gamma_s(V_T, h)\frac{\partial S_a}{\partial V_T}^{-1}$ must be nonnegative to satisfy the necessary condition in Eq. (46) because $-m\gamma/\tilde{D}^2$ is always positive during descent. By this analysis, the possible boundary arcs that satisfy the necessary condition are shown in Fig. 5. The bold lines represent the possible boundary arcs that can be a part of optimal trajectory when the singular arc crosses the state inequality bound. Though the CAS bound case in Eq. (11) is considered only in Fig. 5, the MACH bound case can be analyzed in the same way.

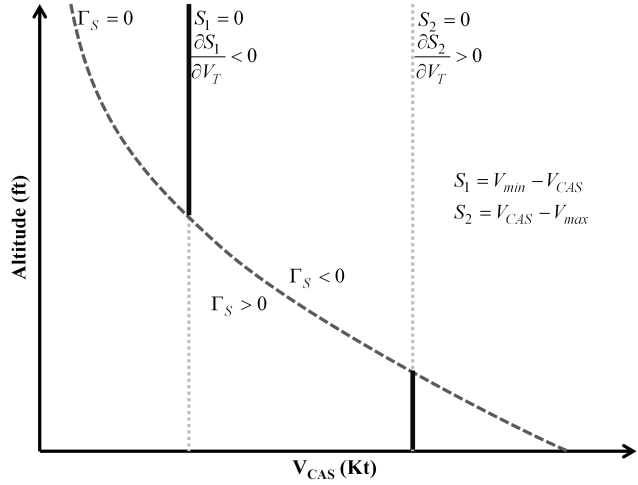

Figure 5. Constrained arc satisfied necessary condition

### 1. Continuity of adjoint variables

The adjoint variables $\lambda_V$ and $\lambda_h$ can be discontinuous for the boundary arc interval $[t_1, t_2]$ from the jump condition in Eq. (35) $\sim$ (37). The adjoint variables can have a discontinuity at any time along the boundary arc interval[23]. From Eq. (33), $\lambda_x$ is constant along the whole optimal trajectory. Furthermore, $\lambda_V$ and $\lambda_h$ are functions of $V_T$ and $h$ in Eq. (43) because $H_\gamma$ and $H_0$ are zero for the boundary arc interval, and hence these are also continuous along the boundary arc because the state variables are continuous. Therefore, the only possible discontinuous points are the junction points. The junction points have two types: the points at which the nonsingular arc and boundary arc are joined and the points at which the singular arc and the boundary arc are joined. To analyze the discontinuity of the adjoint variable at the junction points, we use the following two junction theorems in [26].

**Theorem 2** (Theorem 5.1 in [26]). *Let $t_1$ be the time at which an interior nonsingular arc and a boundary arc of an optimal control $u$ are joined. Let $u^{(r)}$, $r \geq 0$, be the lowest order derivative of $u$ which is discontinuous at $t_1$ and let $p$ be the order of state inequality and $q$ be the order of singular arc. Let $p \leq 2q+r$. If $\nu(t_1) > 0$, then $p + r$ is an even integer.*

**Theorem 3** (Theorem 5.4 in [26])**.** *Let $t_1$ be a point where an interior singular arc and a boundary arc of an optimal control $u$ are joined. Let $q$ be the order of the singular arc and assume that the strengthened GLC-condition holds. Let $u^{(r)}$, $r \geq 0$, be the lowest order derivative of $u$ which is discontinuous at $t_1$ and let $p \leq 2q + r$. Then $\nu(t_1) = 0$, and $q + r$ is an odd integer.*

In all the state inequalities in Eq. (11) and (12), $p = 1$ and the order of singular arc $q = 1$. Therefore $p \leq 2q + r$ holds always in this paper. On the nonsingular arc, $\gamma$ is either $\gamma_{\min}$ or $\gamma_{\max}$, and $\gamma$ stays in the interior of $\Omega$ set on the boundary arc by the assumption. Hence, $\gamma$ at the junction points, at which the nonsingular and the boundary arc are joined, is discontinuous. Therefore, $r = 0$ in this case, and $p + r$ is odd number. According to Theorem 2, if $p+r$ is odd number, then $\nu(t_1) = 0$, which means adjoint variables are continuous at this type of junction points. According to Theorem 3, $\nu = 0$ at the second type of junction points. Furthermore, $r$ should be even number because $q = 1$ in this case. The singular arc in Fig. 3 shows that $r = 0$, and hence this condition is satisfied. By this analysis, we can conclude that adjoint variables are continuous along the optimal trajectory.

## IV. Optimal Trajectory Generation

In this section, we present the optimal trajectory generation algorithm based on the analysis results in the previous section. In this algorithm, we generate the optimal trajectory by combining nonsingular, singular, and boundary arcs. The trajectory is generated by forward and backward integration without solving the optimal control problem. Therefore, it can be implemented in the FMS without additional computational capability.

### A. Algorithm

The proposed optimal trajectory generation algorithm is as follows:

STEP 1. Find the singular arc $(V_T, h)$ using Eq. (44) (For this step, we solved the equation $\Gamma_s(V_T, h) = 0$ numerically).

STEP 2. Find the initial segment of the optimal trajectory by forward integrating the equations of motion in Eq. (1) $\sim$ (3) from the initial point with either $\gamma_{max}$ or $\gamma_{min}$ in $\Omega(V_T, t)$. The control input is determined so that trajectory moves toward the singular arc. The integration stops when the trajectory reaches either the possible boundary arc or the singular arc.

STEP 3. Find the last segment of the optimal trajectory by backward integrating the equations of motion from terminal point with either $\gamma_{max}$ or $\gamma_{min}$ in $\Omega(V_T, t)$. The control input is determined so that trajectory moves from singular arc to the terminal point. The backward integration stops when the trajectory reaches either the possible boundary arc or the singular arc.

STEP 4. Determine the trajectory that is either on the possible boundary arc or the singular arc from section III.B and C. To do this, the trajectory is backward integrated from the junction point obtained from STEP 3 until it reaches the junction point from STEP 2.

In STEP 4, if there exists a junction point between the boundary arc and singular arc, the trajectory changed from the singular arc to the boundary arc or vice versa depending on the possibility.

The expected optimal trajectory structure is shown in Fig. 6. The gray dash lines represent the possible boundary arcs by the optimality conditions in section III. The structure of the optimal trajectory depends on the initial and the terminal conditions as shown in Fig. 6. Depending on the boundary conditions, boundary arc intervals may or may not exist. In case 1 of Fig. 6, both junction points in STEP 2 and STEP 3 are on the boundary arcs, and hence the trajectory generated in STEP 4 consists of three segments: boundary arc, singular arc, boundary arc. On the other hand, the trajectory of STEP 4 is on the singular arc only in case 2 of Fig. 6. The number of segments and the type of segments differ from the boundary conditions. The structure of the optimal solution also depends on the wind profile because the singular arc is varied by the wind condition.

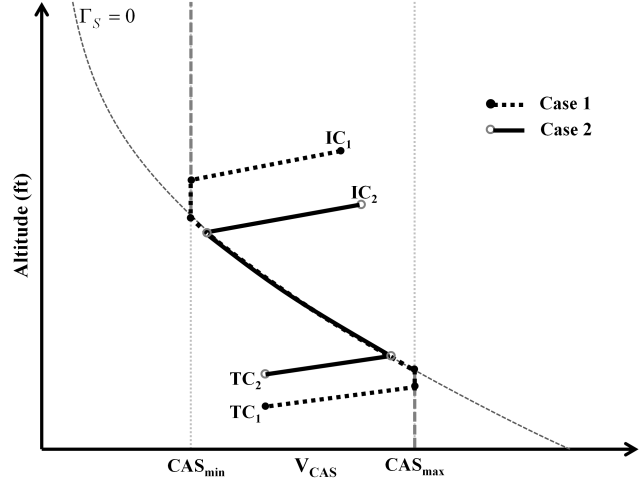

Figure 6. Switching structure of the optimal trajectory

## B. Optimality Check

In this trajectory generation algorithm, we assume that the initial and the last parts of the trajectory, which moves toward the singular arc in STEP 2 and 3, are the parts of the optimal solution. To check if these parts are satisfied the necessary conditions, the sign of $H_\gamma$ should be known. In this paper, $H_\gamma$ is checked by the similar way in [14]. From the necessary condition, the adjoint variables $\lambda_V$ and $\lambda_h$ are continuous at the junction point. Therefore, the values of $\lambda_V$ and $\lambda_h$ at junction point are known from Eq. (28) and (42). From Eq. (34), $\lambda_V$ can be expressed as a function of $V_T$, $h$, and $\lambda_h$. Then Eq. (25) is expressed as follows:

$$\dot{\lambda}_h = K_{cr}\left(\frac{\partial c}{\partial h}V_T + \frac{dW_h}{dh}\right) - \frac{\partial K_{des}}{\partial h} + \frac{\lambda_h V_T \gamma + (K_{des} - K_{cr}(cV_T + W_h))}{\tilde{D} + m\left(g + V_T \mathbf{W}_{h,\chi}\right)\gamma}\left(\frac{\partial \tilde{D}}{\partial h} + m\gamma V_T \frac{\partial \mathbf{W}_{h,\chi}}{\partial h}\right). \tag{49}$$

From the continuity of $\lambda_h$, the values of $\lambda_h$ at the junction points in STEP 2 and STEP 3 are known values. The value at the junction point in STEP 2 is the final value of the initial segment, and the value of STEP 3 is the initial value of the final segment. By integrating Eq. (49) with the known values at the junction points to the backward direction along the initial segment or to the forward direction along the final segment, $\lambda_h(t)$ can be calculated. From the relation between $\lambda_V$ and $\lambda_h$, $\lambda_V$ also can be calculated. Therefore, the time history of $H_\gamma$ is obtained from Eq. (39), and the optimality condition can be checked by Eq. (38). The trajectory generated by the proposed algorithm satisfies this necessary condition, as will be shown by numerical examples.

## V. Numerical Example

To verify that the trajectory generated by the proposed algorithm is the optimal, we compare it to the trajectory generated by the numerical method. For the numerical solution, we use the pseudospectral method since this method provided good solutions in the similar problems in [5]. GPOPS[27, 28] and SNOPT[29] are used as solvers for pseudospectral method and Nonlinear Programming (NLP) problem.

We also analyze the wind effect on the optimal trajectory. The effects of the magnitude of the wind, wind shear, and the wind direction are evaluated. The numerical evaluations of the proposed algorithm are conducted with two aircraft types: B737-500 (B735) and B767-400 (B764). We use BADA[19] performance model for this numerical evaluation. However, since $K_{des}$ and $\tilde{D}$ are functions of $V_T$ and $h$, which are quite general, various performance models can be used to the formulation derived in section III. The boundary condition of this problem is given in Table 1. The maximum range $d_{max}$ is set to -150 NM, which is large enough to satisfy the assumption in Lemma 1. The numerical values of the path constraints in Eq. (11) $\sim$ (14) are given in Table 2.

Table 1. Boundary condition for numerical example

| Initial Condition | | | Final Condition | | |
|---|---|---|---|---|---|
| $V_{CAS}$ (Kt) | $h$ (ft) | $x_s$ (NM) | $V_{CAS}$ (Kt) | $h$ (ft) | $x_s$ (NM) |
| 265 | 35,000 | ¿ $d_{max}$ | 250 | 13,000 | -40 |

Table 2. Path constraints of B735 and B764

| | B737-500 | | B767-400 | |
|---|---|---|---|---|
| | min | max | min | max |
| MACH | 0.45 | 0.82 | 0.45 | 0.84 |
| CAS (Kt) | 220 | 340 | 230 | 360 |
| ROD ($m/s$) | -25.0 | -2.54 | -25.0 | -2.54 |
| $\gamma$ (deg) | -6 | 0 | -6 | 0 |

## A. Fuel Optimal Trajectory

A comparison of the result with the proposed algorithm and with the numerical method is shown in Fig. 7. Zero wind is assumed in this comparison, however, we compared the results for various wind conditions, and the results for all the cases are similar. The two results in Fig. 7 are almost the same except for a very short interval near junction points calculated in STEP 1 and STEP 2 in the proposed algorithm. In the numerical method case, $k\dot{\gamma}^2$ term is added to the Lagrangian cost term in Eq. (19) to prevent the chattering phenomenon that usually occurs when the singular arc exists in the optimal trajectory [30]. This comparison shows that the proposed algorithm and the pseudospectral method generate the same trajectory. We checked the sign of $H_\gamma$ using Eq. (49), and the results are described in Fig. 8. In both initial and the final segments, the control inputs are $\gamma_{max}$ in Fig. 7. Therefore, $H_\gamma$ should be negative to satisfy the necessary condition. The results show that the trajectory from the proposed method satisfies the necessary condition.

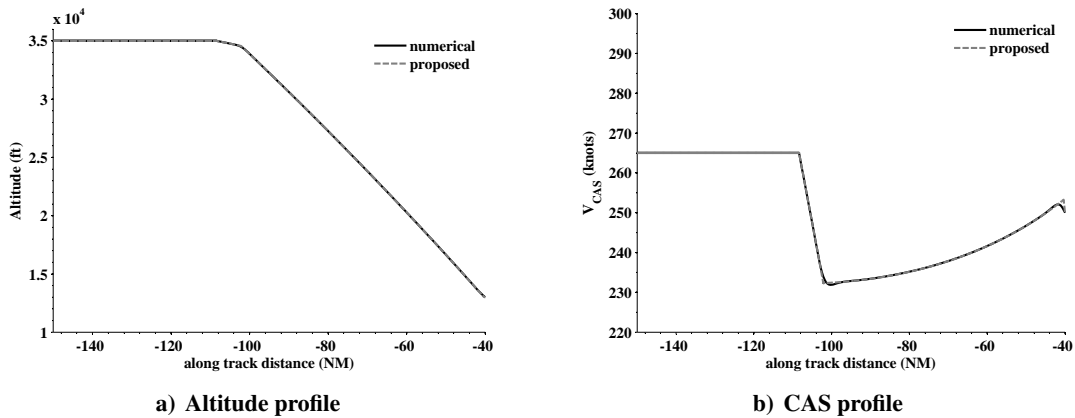

a) Altitude profile      b) CAS profile

Figure 7. Comparison of numerical result and proposed algorithm result (fuel optimal case)

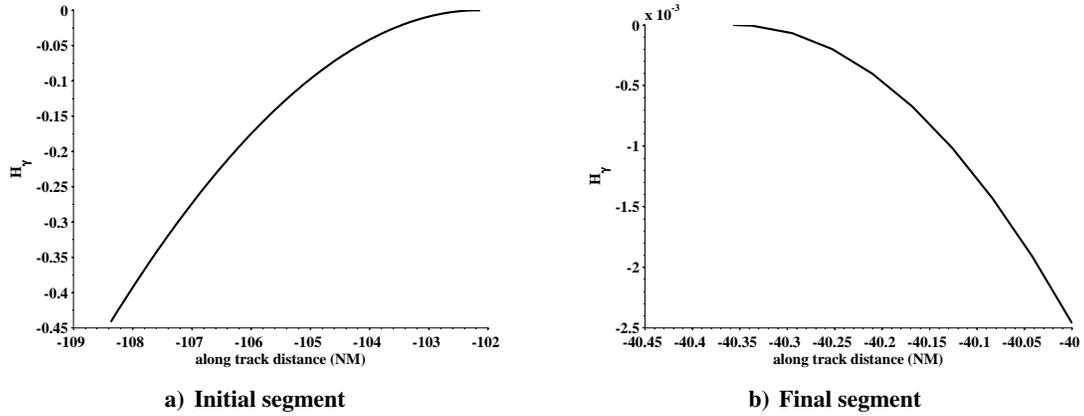

a) Initial segment

b) Final segment

Figure 8. $H_\gamma$ check for STEP 2 and STEP 3 in proposed algorithm

*1. Wind effect*

The fuel optimal trajectories of B735 and B764 for various wind speeds are shown in Fig. 9 and 10. The trajectory performances of the two aircraft are presented in Table 3 and 4. In these examples, wind was assumed to be a constant along track wind to evaluate wind speed effect. In both aircraft, the optimal TOD moves farther from runway threshold as wind speed increases. The descent speed during the singular arc decrease as the tail wind speed increases.

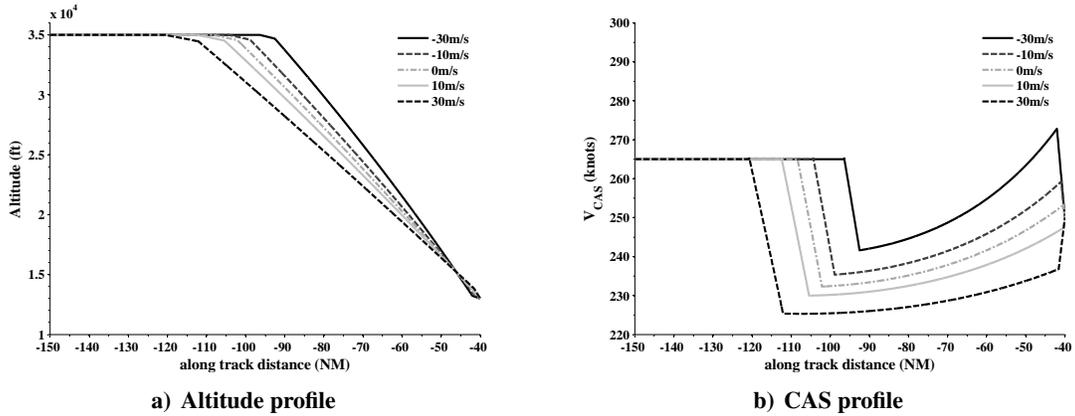

a) Altitude profile

b) CAS profile

Figure 9. B735 fuel optimal trajectories with various wind conditions

Table 3. B735 fuel optimal trajectory performance

| Wind(m/s) | TOD(NM)  | TA(s)    | Fuel(kg) |
|-----------|----------|----------|----------|
| -30       | -96.441  | 1154.489 | 411.314  |
| -10       | -104.336 | 1073.759 | 342.119  |
| 0         | -108.369 | 1038.248 | 311.588  |
| 10        | -112.459 | 1004.226 | 283.161  |
| 30        | -120.711 | 942.495  | 232.326  |

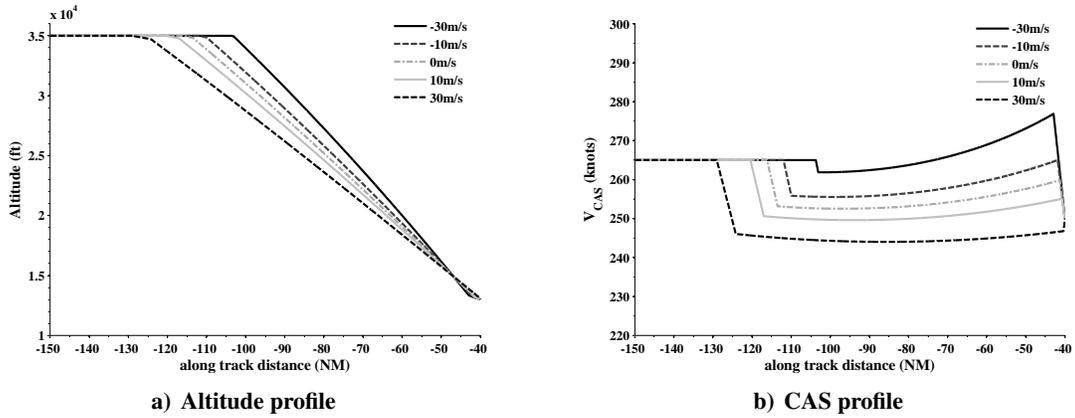

a) Altitude profile  b) CAS profile

Figure 10. B764 fuel optimal trajectories with various wind conditions

Table 4. B764 fuel optimal trajectory performance

| Wind(m/s) | TOD(NM) | TA(s) | Fuel(kg) |
|---|---|---|---|
| -30 | -103.700 | 1139.405 | 729.267 |
| -10 | -111.959 | 1057.429 | 590.346 |
| 0 | -116.167 | 1020.959 | 528.991 |
| 10 | -120.419 | 986.985 | 472.219 |
| 30 | -129.029 | 925.355 | 370.504 |

*2. Wind shear effect*

The fuel optimal CAS profiles of B735 for various wind profiles are shown in Fig. 11. The left column of Fig. 11 is the tail wind case and right column is the head wind case. To evaluate the effects of wind shear, the wind speeds at the cruise altitude are the same for all wind profiles. The influence of the wind shear on the optimal CAS profiles is relatively small. The results show that TOD points are almost the same for the various wind shear cases if the wind speeds at the cruise altitude are the same.

*3. Cross wind effect*

We analyzed the wind direction effect on the performance of the fuel optimal trajectories. For this analysis, we used two constant wind speed, 10 m/s and 30 m/s. The performance difference of B735 and B764 due to the wind directions are shown in Fig. 12. Time of Arrival (TA) and fuel burn increase as the along track wind speed decreases in both B735 and B764 cases.

In order to compare the difference in optimal CAS profile, we compared the two wind profile cases. The first wind profile has the only along track component, and the second wind profile has the same along track wind speed as the first case and has a nonzero cross wind speed. The fuel optimal CAS profile for different cross wind profiles are shown in Fig. 13. In the case where the along track wind speed is zero, the differences in the optimal CAS profile between two cross wind speed cases are small, yet the differences are large in the case where the along track wind speeds are both 30 m/s, but the cross wind components are different (one is zero cross wind, and another has 60 m/s in magnitude, hence the angle $\psi_w = 60 \deg$). In the second comparison case, which Fig. 13(b), the performance differences for the B735 are 15 sec in TA and 15 kg in fuel burn, while for the B764, they are 15 sec in TA and 35 kg in fuel burn. This result shows that we must take the cross wind component into account if we are to determine the true optimal trajectory.

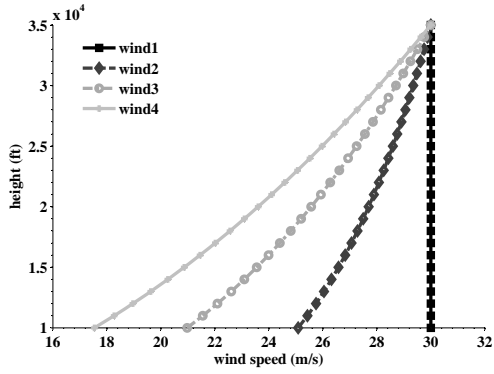

a) Tail wind

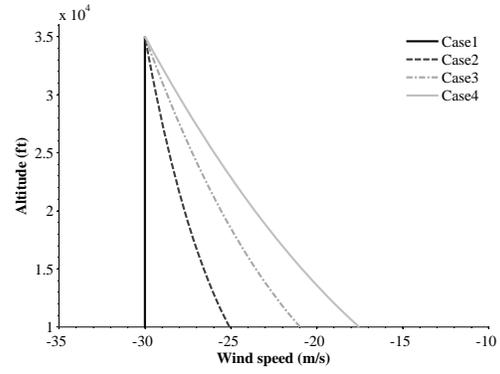

b) Head wind

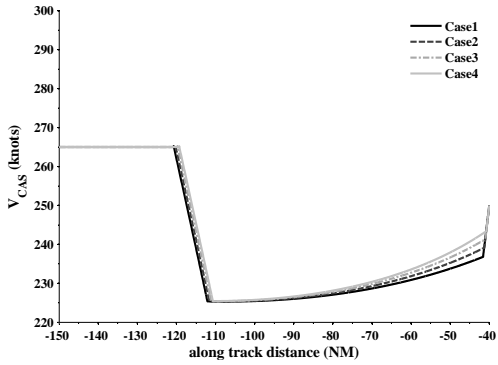

c) B735 T.W

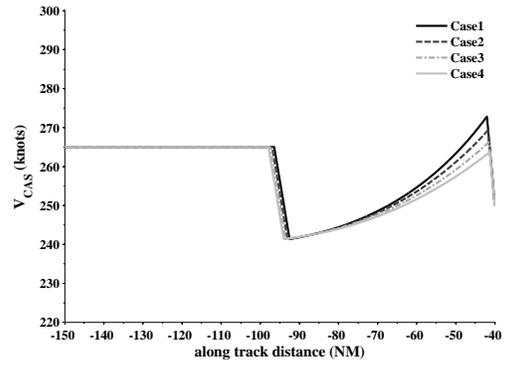

d) B735 H.W

**Figure 11. Wind shear effects: left column - tail wind case, right column - head wind case**

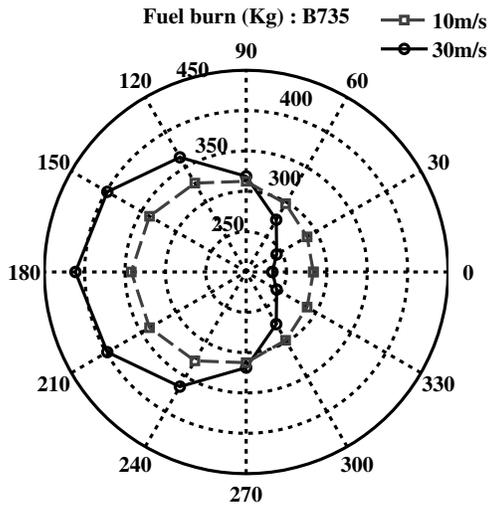

a) B735 Fuel

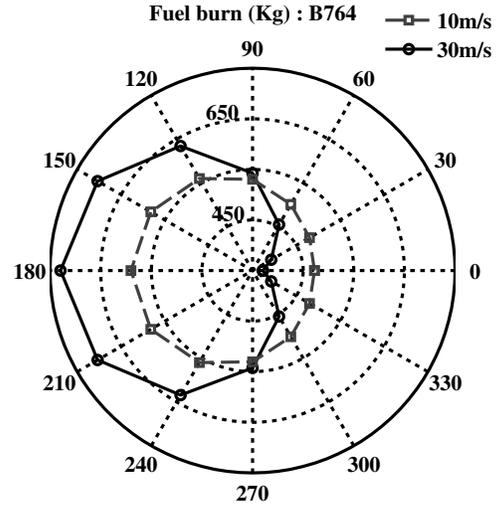

b) B764 Fuel

**Figure 12. Cross wind effect on fuel optimal trajectories**

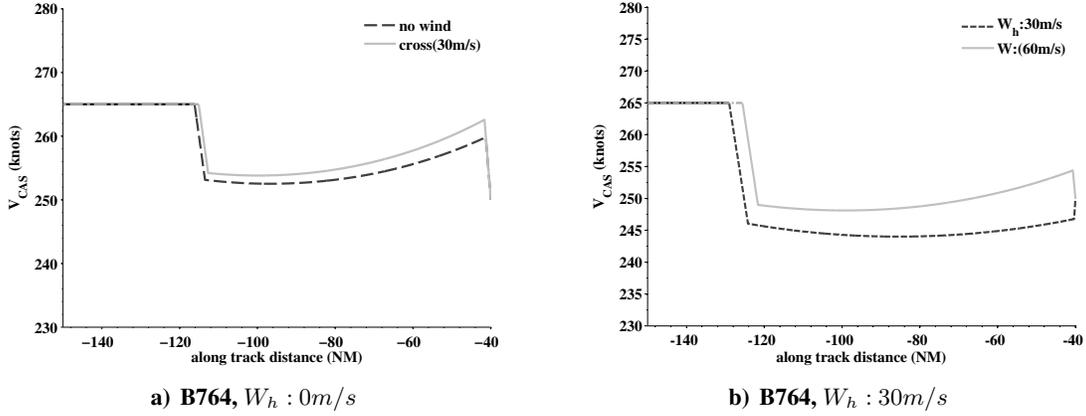

a) **B764**, $W_h : 0m/s$  
b) **B764**, $W_h : 30m/s$

Figure 13. Fuel optimal CAS profiles with same along track winds

## B. Emissions Optimal Trajectory

A comparison of the trajectory from the proposed algorithm and the optimal trajectory from the numerical method for minimizing $NO_x$ emissions is presented in Fig. 14. In this example, wind is assumed to be zero. The two trajectories are almost the same except small deviation at the junction point between nonsingular arc and singular arc. As in the fuel optimal case, the $k\dot{\gamma}^2$ term is added with $k = 0.01$ to protect chattering for the numerical solution. It makes the solution smooth at the junction points compared to the solution from the proposed method. Since the $NO_x$ costs of the two trajectories are the almost the same, and the trajectories are almost the same, this smoothing effect at the junction points is negligible. Hence, it can be concluded that the proposed algorithm and the psuedospectral method generate the same trajectory.

The minimum $NO_x$ trajectories for the B735 and B764 at various wind speeds are shown in Fig. 15 and 16, respectively. The corresponding performance measures are presented in Table 5 and 6. To evaluate wind speed effect, we assume that there is a constant along track wind. The minimum $NO_x$ CAS profiles are different from the fuel optimal speed profiles. The CAS speeds are almost constant along the singular arcs. The TOD points are a little farther than that for the fuel optimal cases, and the arrival times are slightly increased compared to the fuel optimal cases.

We analyzed the wind direction effect on the minimum $NO_x$ trajectory. The test conditions are the same as the minimum fuel trajectory case. The performances of the optimal solutions with the various wind directions are shown in Fig. 17. As maybe seen, the $NO_x$ emission increases as the angle $\psi_w$ in Fig. 2 increases.

The cross wind effect on the minimum $NO_x$ trajectory is shown in Fig. 18. The comparison conditions are the same as the minimum fuel case in Fig. 13(b). The differences in the optimal CAS profile of the two wind profiles are large even though the along track wind speed of both cases are the same as $30m/s$. Furthermore, as shown in Fig. 18(a), the structure of the optimal speed profile can be different depending on the magnitude of the cross wind. The nonzero cross wind case has the boundary arc while the zero cross wind case does not. This result further indicates the importance of considering the cross wind term when deriving the optimal trajectory, regardless of objective.

**Table 5. B735 minimum NOx trajectory performance**

| Wind(m/s) | TOD(NM) | TA(s) | $NO_x$(g) |
|---|---|---|---|
| -30 | -97.008 | 1195.025 | 3698.059 |
| -10 | -104.858 | 1107.863 | 2932.396 |
| 0 | -108.855 | 1068.384 | 2594.677 |
| 10 | -112.890 | 1031.335 | 2282.597 |
| 30 | -121.035 | 960.993 | 1725.159 |

**Table 6. B764 minimum NOx trajectory performance**

| Wind(m/s) | TOD(NM) | TA(s) | $NO_x$(g) |
|---|---|---|---|
| -30 | -104.017 | 1169.034 | 6700.151 |
| -10 | -112.262 | 1083.306 | 5123.663 |
| 0 | -116.457 | 1044.782 | 4428.391 |
| 10 | -120.692 | 1008.701 | 3785.831 |
| 30 | -129.261 | 943.040 | 2636.919 |

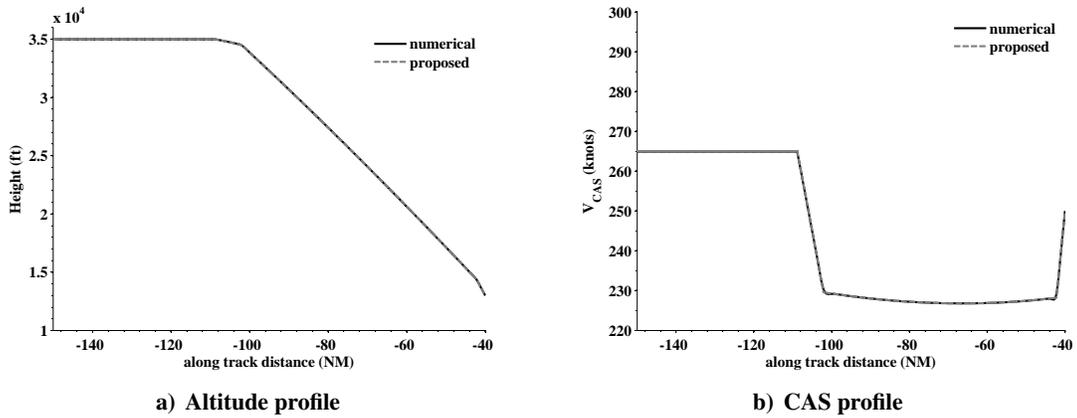

a) Altitude profile  b) CAS profile

**Figure 14. Comparison of numerical method and proposed algorithm results ($NO_x$ optimal case)**

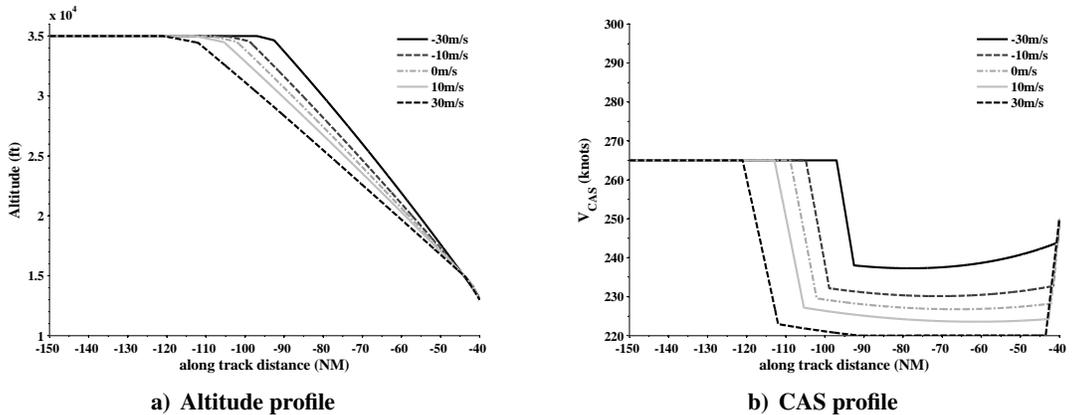

a) Altitude profile  b) CAS profile

**Figure 15. B735 minimum NOx trajectories with various wind conditions**

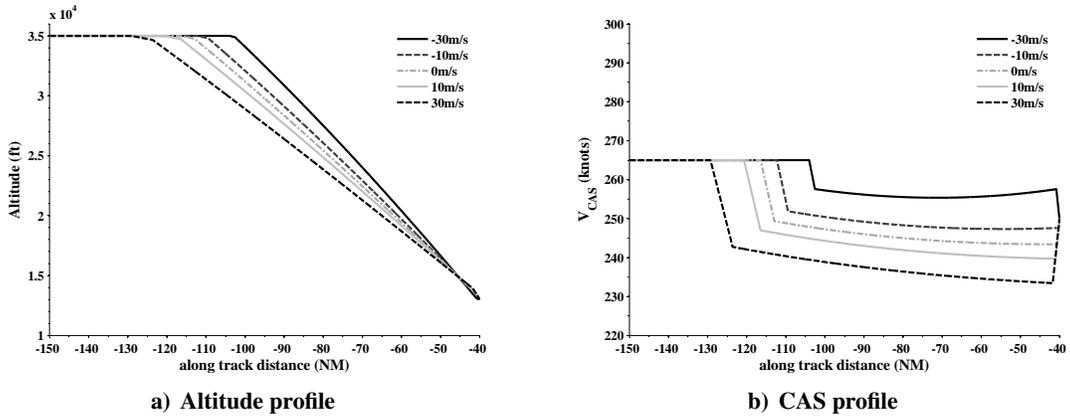

**Figure 16. B764 minimum NOx trajectories with various wind conditions**

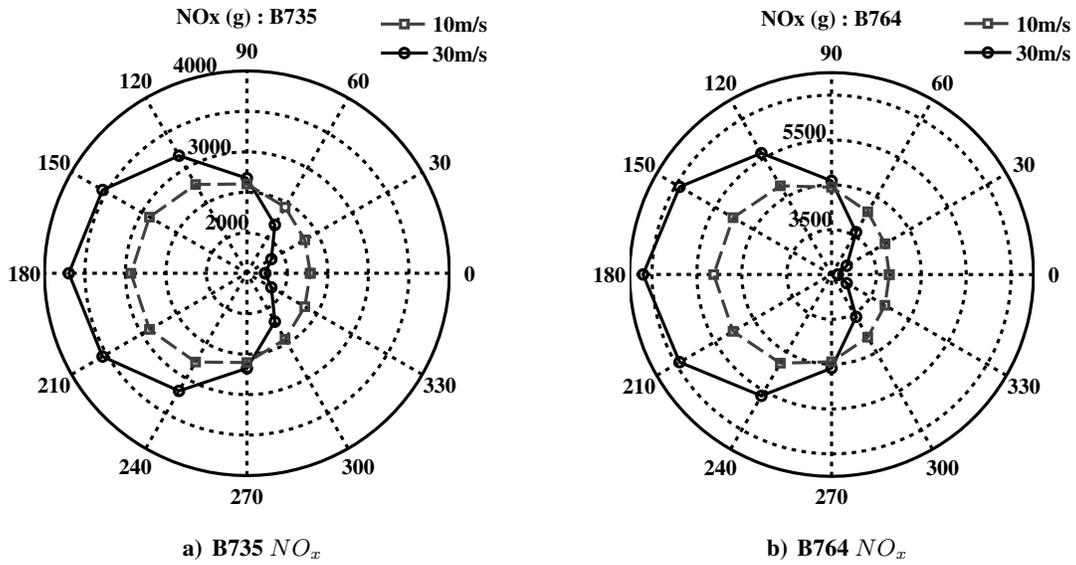

**Figure 17. Cross wind effect on $NO_x$ optimal trajectories**

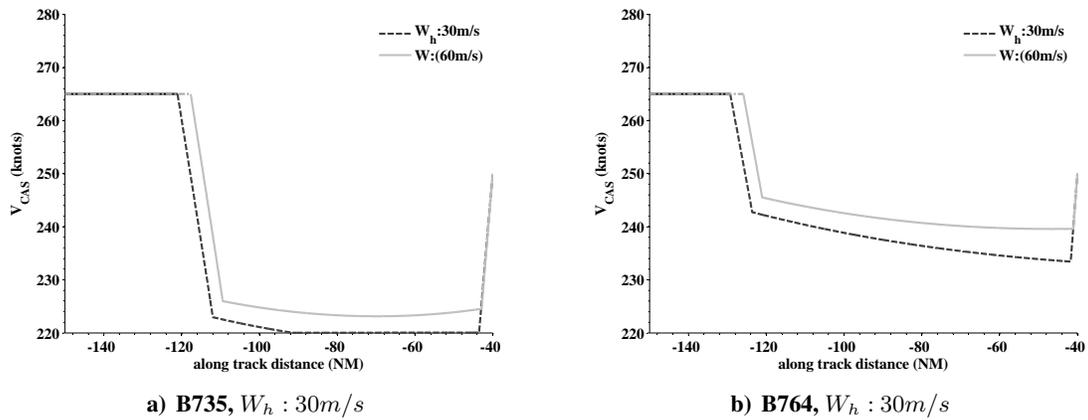

**Figure 18. $NO_x$ optimal CAS profiles with same along track winds**

## VI. Conclusion

We analyzed the optimal vertical trajectories that minimize the environmental impacts of a CDA in the presence of both along track and cross wind components. In this paper, both cross wind and along track wind were assumed to be functions of altitude. Flight idle thrust was assumed during the entire descent phase. The flight range was specified from a point during the latter stages of the cruise segment to the meter fix. We considered CAS and Mach constraints(which are the pure state path constraints), flight path angle constraints, and a maximum descent rate limit(which is a mixed input and state path constraint).

We formulated this problem as a single phase initial condition free optimal control problem. The formulation of the cost functional that we used is quite general since it is a function of airspeed and altitude, and this formulation was used for the fuel cost and the gaseous emissions cost minimization problems. We found an explicit formulation of the singular by using the necessary condition for the optimality of the optimal control problem including path constraints. We also found the existing condition of the boundary arc on which CAS/Mach is constant. From these analyses, we proposed an algorithm to generate the optimal trajectory that minimizes the cost functional by forward and backward integration. The optimal trajectories generated by the proposed algorithm were evaluated by comparing to the optimal trajectories obtained from the numerical method. Since the proposed algorithm does not require additional computational power relative to the current trajectory generation method using the VNAV function in the FMS, this algorithm could easily be implemented in the FMS for the online trajectory optimization.